\documentclass[reqno,12pt]{amsart}
\usepackage{amsmath,amsthm,amssymb,tabularx}

\theoremstyle{plain}

\newtheorem{tm}{Theorem}[section]
\newtheorem{lm}{Lemma}[section]
\newtheorem{cor}{Corollary}[section]
\theoremstyle{definition}
\newtheorem{ex}{Example}[section]

\newcommand{\beq}{\begin{equation}}
\newcommand{\eeq}{\end{equation}}
\newcommand{\bit}{\begin{itemize}}
\newcommand{\eit}{\end{itemize}}
\newcommand{\btm}{\begin{tm}}
\newcommand{\etm}{\end{tm}}
\newcommand{\blm}{\begin{lm}}
\newcommand{\elm}{\end{lm}}
\newcommand{\bcor}{\begin{cor}}
\newcommand{\ecor}{\end{cor}}
\newcommand{\bex}{\begin{ex}}
\newcommand{\eex}{\end{ex}}
\newcommand{\bpr}{\begin{proof}}
\newcommand{\epr}{\end{proof}}

\def\N{\mathbb{N}}
\def\R{\mathbb{R}}
\def\T{\mathbb{T}}
\def\L{{\mathcal L}}
\def\V{{\mathcal V}}
\def\e{\varepsilon}
\let\a\alpha
\def\trans#1#2{{}^{#1}\mkern-3mu #2}

\def \le {\leqslant}
\def \ge {\geqslant}
\let \<\langle
\let \>\rangle
\let\phi\varphi
\let\trg\triangledown
\def\pw{{p,w}}
\def\qw{{q,w^{-1}}}
\def\lpw{\L_p^w(G)}

\def\supp{\operatorname{supp}}
\def\udc#1{\noindent\makebox[\textwidth]{\small UDC #1%
                 \hfill}\vskip10pt}

\begin{document}

\title
{INVARIANT WEIGHTED ALGEBRAS $\L_p^w(G)$}
 \author{Yu. N. Kuznetsova}
 \address{VINITI, Mathematics Department, Usievicha str. 20, Moscow, 125190}
 \email{jkuzn@mccme.ru}
\thanks{Supported by RFBR grant no.~05-01-00982.}

\udc{517.986}
\maketitle

\begin{abstract}
We deal with weighted spaces $\L_p^w(G)$ on a locally compact group $G$.
If $w$ is a positive measurable function on $G$ then we define the space
$\L_p^w(G)$, $p\ge1$, by equality $\L_p^w(G)=\{f:fw\in \L_p(G)\}$. We consider
weights $w$ such that these weighted spaces are algebras with respect to usual
convolution. We show that for $p>1$ such weights exists on any sigma-compact
group. We prove also under minimal requirements a criterion known earlier in
special cases: $\L_1^w(G)$ is an algebra if and only if $w$ is
submultiplicative.
\end{abstract}

Throughout the paper $G$ is a locally compact group, all integrals are taken with respect
to a left Haar measure $\mu$, $p\ge1$, $1/p+1/q=1$ (if $p=1$ then $q=\infty$).
We call any positive measurable function a weight.
Weighted space $\L_p^w(G)$ with the weight $w$ is defined as $\{f:fw\in \L_p(G)\}$,
norm of a function $f$ being $\|f\|_{p,w}=\big(\! \int |fw|^p\, \big)^{1/p}$. Indices $p,w$
are sometimes omitted.

Sufficient conditions on a weight function to define an algebra $\L_p^w(G)$
with respect to usual convolution,
$f*g(s) = \int f(t)g(t^{-1}s)dt$, are well-known.
For $p=1$ it is submultiplicativity:
\beq\label{eq_submult}
w(st) \le w(s)w(t),
\eeq
and for $p>1$ the following inequality (pointwise almost everywhere):
\beq\label{eqw}
w^{-q} * w^{-q} \le w^{-q}.
\eeq
Note that if \eqref{eq_submult} or \eqref{eqw} holds with a constant $C$ (after $\le$ sign)
then for the weight $w_1=Cw$ the same inequality holds without any constant.
Multiplication of a weight by a number changes by the same number the norm of $\L_p^w(G)$,
preserving all the properties of the space. Thus we introduce the notion
of equivalent weights: $w_1$ and $w_2$ are equivalent if with some $C_1$, $C_2$
locally almost everywhere
\beq\label{equiv}
C_1\le {w_1\over w_2} \le C_2.
\eeq

For $p>1$ it is convenient to introduce a dual function $u=w^{-q}$, then the inequality \eqref{eqw}
takes the following form, independent on $p$ and $q$:
\beq\label{equ}
u * u \le u.
\eeq
It is easy to notice that any function $u$ satisfying \eqref{equ} defines a family of
weighted algebras $\L_p^{w_p}(G)$ for all $p\in(1,+\infty)$: $w_p=u^{-1/q}$.

\section{Criterion for the group}

For $p=1$ weighted algebras exist, of course, on any locally compact group
(at least with a unit weight). For $p>1$ we cannot take an arbitrary group, and more
precisely, the following theorem holds:

\btm\label{sigma} For a locally compact group $G$ the following conditions are equivalent:
\bit \item[(i)] $G$ is $\sigma$-compact;

\item[(ii)] for some $p>1$ there exist a weight $w$ satisfying \eqref{eqw}
  (the space $\L_p^w(G)$ is then a convolution algebra);

\item[(iii)] for any $p>1$ there exist a weight $w$ satisfying \eqref{eqw}.
\eit
For an abelian $G$ these conditions are also equivalent to the following:
\bit
\item[(iv)] for some $p>1$ there exist a weight $w$ such that $\L_p^w(G)$ is a convolution algebra;

\item[(v)] for any $p>1$ there exist a weight $w$ such that $\L_p^w(G)$ is a convolution algebra.
\eit
\etm

\bpr Implications (iii)$\Rightarrow$(ii)$\Rightarrow$(iv) and
(iii)$\Rightarrow$(v)$\Rightarrow$(iv) are obvious and do not depend on commutativity of $G$.
We prove that (ii)$\Rightarrow$(i), (i)$\Rightarrow$(iii) and for an abelian group (iv)$\Rightarrow$(i).

(ii)$\Rightarrow$(i). If \eqref{eqw} holds then for some $x$ the integral
$(w^{-q}*w^{-q})(x)=\int w^{-q}(y)w^{-q}(y^{-1}x)dy$ of a strictly positive function is finite.
This implies that $G$ is $\sigma$-compact.

(iv)$\Rightarrow$(i). By \cite[theorem 3]{kuz} there exists an algebra
$\L_p^w(G)$ where $w^{-q}\in \L_1(G)$. Since $w^{-q}$ is positive, $G$ must be $\sigma$-compact.

(i)$\Rightarrow$(iii). We construct a function on $G$ satisfying \eqref{equ}.
Pick a positive function $u_1\in\L_1(G)$ (it exists because $G$ is $\sigma$-compact).
We may assume that $\|u_1\|_1=1$. Define inductively functions $u_n$, $n\in\N$:
$$
u_{n+1} = u_1*u_n.
$$
Clearly $\|u_n\|_1\le1$ for all $n$.
We put now $u=\sum n^{-2} u_n$ and prove that \eqref{equ} holds.
Note the following elementary fact:
\begin{gather*}
\sum_{n=1}^{m-1} {1\over n^2(m-n)^2} \le 8\zeta(2){1\over m^2}.
\end{gather*}
Estimate now the convolution $u*u$:
\begin{gather*}
u*u
 = \sum_{n,k=1}^\infty {u_n*u_k\over n^2\,k^2}
 = \sum_{n,k=1}^\infty {u_{n+k}\over n^2\,k^2} = \\
=\sum_{m=1}^\infty \sum_{n=1}^{m-1} {u_m\over n^2(m-n)^2}
\le 8\zeta(2) \sum_{m=1}^\infty {u_m\over m^2} \equiv Cu.
\end{gather*}
Changing $u$ to $u/C$, we get \eqref{equ}.
\epr

\section{Technical lemmas}

In general any positive measurable function can be taken as a weight.
Some authors assume also that $w^p$ is locally summable, but this is redundant if
$\L_p^w(G)$ is an algebra:

\blm\label{int}
If the space $\L_p^w(G)$, $p\ge 1$, is a convolution algebra then $w^p$ is locally summable.
\elm

\bpr
Take a compact set $A\subset G$ of positive measure. Here and further $I_A$ denotes
the characteristic function of a set $A$.
Consider the functions $\phi=I_A/\max\{1,w\}$ and $\psi=I_{A^{-1}\cdot A}/\max\{1,w\}$.
As $\phi,\psi\in \L_p(w)$, then also $\tau=\phi*\psi\in \L_p(w)$.
At the same time $\phi\in \L_1$ and $\psi\in \L_\infty$,
so that $\tau$ is continuous. Since also $\tau|_A>0$ we have $\min_A \tau=\tau_0>0$.
Thus, $I_A\le \tau/\tau_0$, what implies $I_A\in \L_p(w)$ or, equivalently, $w\in \L_p(A)$.
\epr

\blm\label{compact_dense}
The set $B_0(G)$ of all bounded compactly supported functions
is a dense subspace of every algebra $\lpw$.
\elm

\bpr
Inclusion $B_0(G)\subset\lpw$ follows from the previous lemma.
It is known that $B_0(G)$ is dense in $\L_p(G)$, therefore $w^{-1}B_0(G)$
is dense in $\lpw=w^{-1}\L_p(G)$. Let now $f\in w^{-1}B_0(G)$ be compactly supported,
but not necessarily bounded. Changing it on a set of arbitrary small measure $\e$,
we can make $f$ bounded; $\e$ is to be chosen according to continuity of the
integral $\|f\|_\pw=\|fw\|_p$ as a set function.
\epr

The property of submultiplicativity \eqref{eq_submult} is essential in the weighted algebras
theory because exactly this property posess the weights of $L_1^w(G)$ algebras
(see theorem \ref{submult} below). We need the following lemma on submultiplicative functions
(it is in fact proposition 1.16 of \cite{edw}):

\blm\label{cbound}
Let a measurable function $L:G\to\R$ be submultiplicative, i.e.~satisfy
\eqref{eq_submult}, and positive. Then $L$ is bounded and bounded away from zero on
any compact set.
\elm

The following condition studied first by R.~Edwards \cite{edw}
is also important for the weighted spaces theory.
A weight $w$ is of moderate growth if for all $s\in G$
\beq\label{ledw}
L_s = {\rm ess}\sup_t {w(st)\over w(t)} <\infty.
\eeq
This condition is equivalent to the space $L_p^w(G)$ (for any $p\ge1$) being
translation-invariant \cite[1.13]{edw}. In the non-commutative case, \eqref{ledw}
corresponds to left translations; taking $w(ts)$ indtead of $w(st)$, we get a condition for
right translations, in general not equivalent to the former.
Immediate calculations show that
\beq\label{trans_norm}
\sup_{f\ne0} {\|\trans sf\|_{p,w}\over \|f\|_{p,w}} = L_{s^{-1}}.
\eeq
The condition \eqref{ledw} implies that $L_s>0$, $L_{st}\le L_sL_t$, and
$$
{\rm ess}\inf_t {w(st)\over w(t)} = 1/L_{s^{-1}} >0.
$$

\blm\label{edw_a.e.}
If \eqref{ledw} holds for locally almost all $s\in G$ then it holds for all $s\in G$.
\elm

\bpr
Let $S\subset G$ be the set of $s$ for which the inequality
\eqref{ledw} holds. By assumption $S$ and hence $S^{-1}$ is locally of full measure.
Pick a set $T\subset S\cap S^{-1}$ of positive finite measure.
Then $T\cdot T^{-1}$ contains a neighborhood of identity $U$.
As $L$ (finite or infinite) is submultiplicative, $S$ is closed under multiplication and therefore
$U\subset T\cdot T^{-1}\subset S\cdot S\subset S$. By the same reason $SU\subset S$,
and since $S$ (being locally of full measure) is everywhere dense, then $S=G$.
\epr

\blm\cite[th. 2.7]{feicht}\label{contin}
If a weight $w$ satisfies \eqref{ledw} and is locally summable then it is equivalent to
a continuous function.
\elm

\bcor
Let $L_p^w(G)$ be an algebra with a weight $w$ satisfying \eqref{ledw}.
Then $w$ is equivalent to a continuous function.
\ecor

\bpr
By lemma \ref{int} for any compact set $F$ we have $w\in \L_p(F)\subset \L_1(F)$,
therefore we can apply lemma \ref{contin}.
\epr

On a compact group any continuous function is equivalent to a constant function,
thus on a compact group all translation-invariant weighted algebras are isomorphic
to the usual algebra $\L_p(G)$.
Converse of the corollary does not hold:

\bex\label{not-invar-alg}
There exist an algebra $L_2^w(\R)$ such that $w$ is continuous but does not
satisfy \eqref{ledw}. Let $0<\a_n<1$, $A_n=[n+\a_n,n+1]$.
We put $w|_{A_n}=1+n^2$, $w(n+\a_n/2)=1+|n|$ and extend $w$ piecewise linearly.
For $\a_n=n^{-2}$ the condition \eqref{eqw} is satisfied but \eqref{ledw}
does not hold in any neighborhood of zero.
\eex

\section{Criterion for the algebra $\L_1^w(G)$}\label{sect_l1}

In the case when $p=1$ the class of weights defining convolution algebras $\L_p^w(G)$
admits a complete description, and it turns out that every weight is
equivalent to a continuous function.
The following theorem was proved by Grabiner \cite{grab} in the case
of the real half-line (without statement of continuity which is false on the half-line).
Edwards \cite{edw} proved equivalence of (i) and (ii)
on a locally compact group under assumption of upper-semicontinuity of $w$,
and later Feichtinger \cite{feicht} for translation-invariant algebras $\L_1^w(G)$.
Our theorem generalizes these results.

\btm\label{l1}
For a weight $w$ the following conditions are equivalent:
\bit
\item[(i)] $w$ is equivalent {\rm(}in the sense of \eqref{equiv}{\rm)}
to a continuous submultiplicative function;

\item[(ii)] $\L_1^w(G)$ is a convolution algebra;

\item[(iii)] for some $p$, $1\le p<\infty$, the inclusion
$\L_1^w(G)*\L_p^w(G)\subset \L_p^w(G)$ holds.
\eit
\etm

\bpr
Implications (i)$\Rightarrow$(ii) and (i)$\Rightarrow$(iii) are immediate
whereas (ii)$\Rightarrow$(i) is a special case of (iii)$\Rightarrow$(i) with $p=1$.
We prove therefore only (iii)$\Rightarrow$(i).

Inclusion (iii) implies (cf.\ \cite[38.27]{HR}) that with some constant $C$
$$
\|f*g\|_{p,w} \le C\|f\|_{1,w}\|g\|_{p,w}.
$$
Repeating argument of lemma \ref{int} we conclude that $w^p$ together with $w$
are locally summable. Thus the spaces $\L_1^w(G)$, $\L_p^w(G)$ contain charasteristic
functions of all sets of finite measure. For such sets $A$, $B$ and arbitrary $s,t$
we have pointwise
\beq\label{eqI}
\mu(A) I_{stB} \le I_{sA}*I_{A^{-\!1}tB},
\eeq
whence
\beq\label{x}
\mu(A) \|I_{stB}\|_{p,w} \le C\|I_{sA}\|_{1,w}\|I_{A^{-\!1}tB}\|_{p,w}.
\eeq

We need here a generalization of the Lebesgue differentiation theorem.
On a locally compact group one may state the theorem as follows
(see a general statement in the review \cite{bruck} and specifications for the group case
in \cite{ion}): there exists a family $\V$ of sets of positive measure directed by downward
inclusion such that for any locally summable function $f$
\beq\label{lebe}
\lim_{V\in\V} {1\over \mu(V)} \int_{xV} f(t)dt = f(x)
\eeq
for locally almost all $x\in G$.
At that every $V\in \V$ contains the identity, and every neighborhood of identity contains
eventually all $V\in\V$ \cite[VIII, 1-2]{ion}.

So, for locally almost all $s\in G$ \eqref{lebe} holds with $f=w$, $x=s$.
For each such $s$ \eqref{lebe} holds both with $f=w^p$ and $f=\trans sw^p$
for locally almost all $x=t$. For such $s,t$ and any $\e>0$
for sufficiently small $V\in\V$
\begin{gather*}
\|I_{stV}\|_{p,w}^p=\int_{tV}\big(\trans sw(r)\big)^pdr > w^p(st)\mu(V)/(1+\e),  \\
\|I_{tV}\|_{p,w}^p=\int_{tV}w^p(r)dr < w^p(t)\mu(V)(1+\e).
\end{gather*}
Fix $V$ such that these inequalities hold. Since the integral of
$\trans sw^p$ are continuous as functions of a set,
there is a compact set $B\subset V$ such that
\begin{gather*}
\|I_{stV}\|_{p,w}< (1+\e)\,\|I_{stB}\|_{p,w},  \qquad
\mu(V)<(1+\e) \mu(B).
\end{gather*}
Moreover, there exists a neighborhood of identity (with compact closure) $V_0$ such that
$$
\|I_{V_0^{-1}tB}\|_{p,w}<(1+\e)\|I_{tB}\|_{p,w}.
$$
And, finally, for sufficiently small $A\in\V$, $A\subset V_0$ holds
$$
\|I_{sA}\|_{1,w} <(1+\e)\mu(A)w(s).
$$
Obviously, $\|I_{A^{-1}tB}\|_{p,w}\le \|I_{V_0^{-1}tB}\|_{p,w}$ and
$\|I_{tB}\|_{p,w}\le\|I_{tV}\|_{p,w}$.
Uniting all these inequalities with \eqref{x}, we get:
$$
\mu(A) \mu(B)^{1/p} w(st) < C (1+\e)^{3/\!p\,+3} \mu(A)\mu(B)^{1/p}w(s)w(t),
$$
and in the limit as $\e\to0$
\beq\label{y}
w(st) \le C w(s)w(t).
\eeq
This inequality is obtained for locally almost all $t$ with fixed $s$ for locally almost all $s$.
But by lemma \ref{edw_a.e.} for $w$ the condition \eqref{ledw} holds, and by lemma \ref{contin}
$w$ is equivalent to a continuous function $w_1$. For $w_1$ \eqref{y}
(with another constant) holds for all $t$ and $s$. Finally, multiplying $w_1$ by this constant,
we get a continuous submultiplicative weight.
\epr

As a corollary we obtain a description of multipliers of the algebra $\L_1^w(G)$.
Gaudry \cite{gau} proved under assumption of upper-semicontinuity that multipliers of
$\L_1^w(G)$ may be identified with the weighted space ${\mathcal M}^w(G)$ of regular Borel measures
such that $\int wd\/|\mu|<\infty$.
As the weight can be always chosen continuous, statement of the theorem is simplified:

\btm
A bounded operator $T$ on an algebra $\L_1^w(G)$ commutes with right translations
if and only if it is a convolution with a measure $\mu\in {\mathcal M}^w(G)$: \
$Tf=\mu*f$ for all $f\in \L_1^w(G)$.
\etm

If $p>1$, the weight of an algebra cannot in general be chosen continuous:

\bex\label{countex}
$G=\T$, $p=2$. Parametrize the circle by $t\in {[-1,1]}$ and take $w(t)=|t|^{-1/4}$
(another example: $w(t)=|t|^{1/4}$). $\L_2^{w}(\T)$ is an algebra but the weight
is not equivalent to a continuous function. At that $L_2^w(\T)$ is not invariant under
translations and $L_1^w(\T)$ with the same weight is not an algebra.
\eex

The fact that $L_1^w(G)$ is an algebra implies the weight is submultiplicative and
all the spaces $L_p^w(G)$, $p\ge1$, are translation-invariant (i.e.\ \eqref{ledw} holds).
Converse is true for abelian groups:

\btm\label{submult}
Let $G$ be an abelian locally compact group, $\L_p^w(G)$ an algebra with a weight $w$
satisfying \eqref{ledw}. Then $w$ is equivalent to a submultiplicative function.
\etm

\bpr
By lemma \ref{contin} we may assume $w$ is continuous, and by lemma \ref{cbound}
the function $L$ is bounded on any compact set. Pick a compact set $D=D^{-1}$
of positive measure, and let $L(r)\le N$ for $r\in D$.
Then for $s\in G$, $r\in D$ we have
\beq\label{wD}
{w(s)\over N}\le w(sr) \le Nw(s).
\eeq
Take now arbitrary $s,t\in G$ and use inequality \eqref{eqI} with $V=U=D$.
We get then
$$
\mu(D) \|I_{stD}\| \le \|I_{sD}\|\cdot\|I_{tD^2}\|,
$$
whence by \eqref{wD}
$$
\mu(D)^{1+1/p} w(st)/N \le N^2 \mu(D)^{1/p}\mu(D^2)^{1/p} w(s)w(t),
$$
i.e. $w(st) \le C_1 w(s)w(t)$, what completes the proof.
\epr

Note that on a discrete group weight of any algebra for all $p\ge1$
is submultiplicative. Inequality \eqref{eq_submult} obtains when passing from
$I_{st}=I_s*I_t$ to the norms in $\L_p^w(G)$.

\section{Approximate units}\label{sect_approx}

Algebras $\L_p^w(G)$, as well as those without weight, have a unit iff $G$
is discrete. In the classical case (of a compact group) $\L_p(G)$ always have
approximate units. In the weighted case moderate growth of the weight
(see th. \ref{approx}) is sufficient for an algebra to have an a.u.
If the weight is not moderate, this theorem may not hold, see example \ref{no_a.e.}.
Invariant algebras $\L_p^w(G)$ do not have bounded approximate units (theorem \ref{bae}).
It follows these algebras are not amenable \cite{runde}.
In this section, theorems are proved for left a.u., but the same is true
for right a.u. with interchange of $s$ and $t$ in the condition \eqref{ledw}.

\blm\label{transl}
Let $\L_p^w(G)$, $p\ge1$ be an invariant algebra.
Then for any $f\in \L_p^w(G)$ and $\e>0$ there exists a neighborhood of
identity $U$ such that for all $t\in U$
$\|f-\trans tf\|_{p,w}<\e$.
\elm

\bpr
Suppose first that $f\in\L_p^w(G)$ is compactly supported, $F=\supp f$.
Pick a relatively compact neighborhood of identity $U$, then
$\supp \trans tf\subset UF$ when $t\in U$.
By lemma \ref{contin} we may assume $w$ is continuous, thus bounded on every
compact set, so that $C=\sup\limits_{UF} w<\infty$.
Now
$$
\|f-\trans tf\|_{p,w} = \Big(\int_{UF} |f(x)-f(tx)|^pw^p(x)dx\Big)^{1/p}
 \le C \|f-\trans tf\|_p,
$$
where the latter norm is less than $\e$ for $t$ in a sufficiently small 
neighborhood of identity $V\subset U$.

Let now $f\in\L_p^w(G)$ be arbitrary. For any $\e>0$ it may be approximated by
a compactly supported function $\phi\in\L_p^w(G)$, $\|f-\phi\,\|_{p,w}<\e$ (lemma \ref{compact_dense}).
Pick again a symmetric relatively compact neighborhood of identity $U$.
By lemma \ref{cbound} the function $L$ defined in the formula \eqref{ledw}
is bounded on every compact set, so that $D=\sup_U L<\infty$.
Now by equality \eqref{trans_norm} for $t\in U$
$$
\|\trans tf-\trans t\phi\|_\pw = \|\trans t(f-\phi)\|_\pw\le L_{t^{-1}} \|f-\phi\|_\pw \le D \|f-\phi\|_\pw,
$$
and
$$
\|f-\trans tf\|_\pw \le \|f-\phi\|_\pw + \|\phi-\trans t\phi\|_\pw +
 \|\trans tf-\trans t\phi\|_\pw < \e+D\e+C\|\phi-\trans t\phi\|_p,
$$
what is less than $\e(1+C+D)$ for $t$ in a sufficiently small 
neighborhood of identity $V\subset U$.
\epr

\btm\label{approx}
Let $\L_p^w(G)$, $p\ge1$ be an invariant algebra.
The net $\xi_\nu=I_\nu/\mu(\nu)$, where $\nu$ runs over the net of
all relatively compact neighborhoods of identity, is a left approximate
unit in $\L_p^w(G)$.
\etm

\bpr
Note first that all $\xi_\nu$ belong to $\L_p^w(G)$ (prop. \ref{compact_dense}).
Convergence $\xi_\nu*f\to f$ for every $f\in \L_p^w(G)$ is proved in a
standard way. By lemma \ref{transl} exists a neighborhoods of identity $U$ such
that $\|f-\trans {t^{-\!1}\!}f\|_\pw<\e$ for $t\in U$.
We can estimate the norm $\|\xi_\nu*f-f\|_\pw$, $\nu\subset U$
using functionals of the conjugate space: for every $\phi\in \L_q^{w^{-1}}(G)$
\begin{gather*}
|\langle \xi_\nu*f-f, \phi\rangle| =
\Big|\int_G \big[ (\xi_\nu*f)(x)-f(x)\big] \phi(x) dx\Big|
\le\\
 \le \int_G \int_\nu {|f(t^{-1}x)dt - f(x)|\over \mu(\nu)}dt \,|\phi(x)|dx
 \le {1\over \mu(\nu)} \int_\nu \langle |\trans{t^{-1}}f - f|,|\phi|\rangle dt
\le\\
\le \sup_{t\in \,U}\|\trans {t^{-\!1}\!}f-f\|_{p,w}\|\phi\|_\qw
<\e \|\phi\|_\qw,
\end{gather*}
i.e. $\|\trans {t^{-\!1}\!}f-f\|_{p,w}<\e$, what proves the theorem.
\epr

\bex\label{no_a.e.}
If the weight fails to satisfy \eqref{ledw}, the statement of theorem \ref{approx}
may not be true. The following algebra is a counterexample.
Take $G=\R$, $p=2$. We denote $\bar n=\max(|n|,1)$ and define the weight
in the following way:
$$
w(t)=\begin{cases}
 \bar n,& t\in[\,n,n+1/\bar n^2), \\
 \bar n^2,& t\in[\,n+1/\bar n^2,n+1).
 \end{cases}
$$
We show first that (2) holds for $w$ after a multiplication by some
constant, i.e. $\L_2^w(\R)$ is an algebra. Denote
$$
I_n=I_{[n,n+1/\bar n^2)},\quad I'_n=I_{[n+1/\bar n^2,n+1)}.
$$
In these notations $w=\sum(\bar n I_n + \bar n^2I'_n)$. Using a trivial
estimate $I_A*I_B\le \min\{\mu(A),\mu(B)\} I_{A+B}$ and inequality
$$
\sum_{n=-\infty}^\infty {1\over \bar n^\a(\overline{m-n})^\a} \le
 {2^{\a+1}\over \bar m^\a}(2\sum_{n=1}^\infty {1\over \bar n^\a}+1) = {C_\a\over\bar m^\a}
$$
for $\a=2,4$ and integer $m$, we can estimate convolution in (2):
\begin{gather*}
w^{-2}*w^{-2}=\sum_{n,m} \Big( {1\over \bar n^2\bar m^2} I_n*I_m +
 2 {1\over \bar n^2\bar m^4} I_n*I'_m + {1\over \bar n^4\bar m^4} I'_n*I'_m \Big)
 \le\\\le
\sum_{n,m} \Big(
 {1\over \bar n^2\bar m^2\max(\bar n^2,\bar m^2)} I_{n+m+[0,2)}+
 2 {1\over \bar n^2\bar m^4\bar n^2} I_{n+m+[0,2)}
 + {1\over \bar n^4\bar m^4} I_{n+m+[0,2)} \Big)
 \le\\\le
\sum_k I_{k+[0,2)} \Big( \sum_n {3\over \bar n^4(\overline{{k-n}})^4} +
 {4\over \bar k^2}\sum_n {1\over \bar n^2(\overline{{k-n}})^2}\Big)
 \le\\\le
\sum_k I_{k+[0,2)} \Big( {3C_4 \over \bar k^4} +
 {4\over \bar k^2}\cdot {C_2\over\bar k^2}\Big)
 \le C\sum {1\over\bar k^4}I_{[k,k+1)} \le Cw^{-2}.
\end{gather*}
Obviously $\L_2^w(\R)$ is not translation invariant.

Now we prove that this algebra has no a.u. consisting of nonnegative
functions. Suppose the opposite, i.e. that $e_\a\ge0$ are an a.u.

First we show that $i_\a=\int_{-1/4}^{1/4}e_\a\not\to0$.
Let $I_\delta=I_{[0,\delta]}$ be the indicator function of $[0,\delta]$,
$0<\delta<1/4$. As
$$
(I_\delta*e_\a)(t)=\int_{t-\delta}^t e_\a,
$$
then $I_\delta*e_\a\le i_\a$ if $t\in[0,\delta]$. But if $i_\a\to0$, then
$$
\|I_\delta-I_\delta*e_\a\|_{2,w}^2\ge \int_0^{\delta} (I_\delta-I_\delta*e_\a)^2
 \ge \delta\cdot (1-i_\a)^2\to\delta\ne0.
$$
Thus, $i_\a\not\to0$, what means that integral of $e_\a$ either over
$[-1/4,0]$ or over $[0,1/4]$ does not tend to zero. Suppose the latter
(otherwise we should define $\tilde I_n$ below as left shifts of $I_n$
instead of right ones).

Introduce functions
$$
f=\sum_{n=1}^\infty \a_n I_{2^n},\quad
g=\sum_{n=1}^\infty \gamma_n \tilde I_{2^n},
$$
where $\tilde I_n = I_{n+[1/\bar n^2,2/\bar n^2)}$.
We will choose $\a_n$, $\gamma_n$ so that
$f\in \L_2^w(\R)$, $g\in\L_2^{w^{-1}}(\R)$.
According to the definition of weight
\begin{gather*}
\|f\|_{2,w}^2 = \sum_{n=1}^\infty \a_n^2 2^{2n} {1\over 2^{2n}} = \sum \a_n^2, \\
\|g\|_{2,w^{-1}}^2 = \sum_{n=1}^\infty \gamma_n^2 2^{-4n} {1\over 2^{2n}} =
 \sum \gamma_n^2 2^{-6n}.
\end{gather*}
Thus we can put $\gamma_n=2^{3n}\beta_n$ and take any sequences
$\a,\beta\in\ell_2$.

Now we show that $f*e_\a\not\to f$. It is sufficient to show that
$\<f*e_\a,g\>\not\to \<f,g\>$. Since $\<f,g\>=0$, we should show that
$$
\<f*e_\a,g\> = \<e_\a,f^\trg*g\>\not\to0
$$
(here $f^\trg(t)=f(-t)$).
Let us calculate the convolution
$$
f^\trg*g = \sum \a_n\gamma_k I_{2^n}^\trg*\tilde I_{2^k}.
$$
For fixed $n,k$
$$
\supp I_{2^n}^\trg*\tilde I_{2^k} = \supp I_{2^n}^\trg+\supp \tilde I_{2^k}
 = 2^k-2^n+[-2^{-2n},2^{-2k}].
$$
We will be interested below in the segment $[0,1]$ only, i.e. convolutions
with $n=k$. These we calculate explicitly:
$$
I^\trg_{2^n}*\tilde I_{2^n} = 2^{-2n}J_n,
$$
where
$$
J_n = \min(1+2^{2n}t,2) - \max(2^{2n}t,1) = \begin{cases}
 2^{2n}t, & t\in[0,2^{-2n}],\\
 2-2^{2n}t, & t\in[2^{-2n},2^{1-2n}].\end{cases}
$$
This function is piecewise linear, and $J_n(0)=J_n(2^{1-2n})=0$, $J_n(2^{-2n})=1$.
Thus,
$$
(f^\trg*g)|_{[0,1]} = \sum \a_n \gamma_n 2^{-2n} J_n = \sum \a_n \beta_n 2^n J_n.
$$
Put now $\a_n=\beta_n=2^{-n/2}$. As required, $\a,\beta\in\ell_2$,
and
$$
(f^\trg*g)|_{[0,1]} = \sum J_n = J.
$$
Next we estimate $J$:
\begin{gather*}
J(2^{-2n}) = \sum_{k=1}^n 2^{2k}\cdot 2^{-2n} = {2^{2n}-1\over3}2^{2-2n}
 = {4\over3}(1-2^{-2n}),\\
J(2^{1-2n}) = \sum_{k=1}^{n-1} 2^{2k}\cdot 2^{1-2n} = {2^{2n-2}-1\over3} 2^{3-2n}
 = {2\over3}(1-2^{1-2n})
\end{gather*}
(in the latter case $n>1$).
Remember that all summands are piecewise linear, therefore
$1/3 \le J\le 4/3$ on $(0,1/4]$. It follows that
$$
{1\over3}\int_0^{1/4} e_\a \le \int e_\a J \le {4\over3}\int_0^{1/4} e_\a,
$$
i.e. $\int e_\a J\not\to0$.
Since $\<e_\a,f^\trg*g\> \ge \int e_\a J$, we proved that $e_\a$ is not an
approximate unit.
\eex

\btm\label{bae}
Let $G$ be a non-discrete group, $p>1$, and let $\L_p^w(G)$ be an invariant
algebra. Then this algebra has no bounded approximate unit.
\etm

\bpr Suppose that $\L_p^w(G)$ has a left b.a.u. Then \cite[th. 32.22]{HR}
$\L_p^w(G)\circ X$ is a closed linear subspace in $X$ for any left module $X$
over $\L_p^w(G)$ with a multiplication $\circ$.
Take $X=\L_q^{w^{-1}}(G)$. This is a left module over $\L_p^w(G)$ with
multiplication $f\circ g = g*f^\trg$, $f^\trg(x)=f(x^{-1})$.

Show that $Y=\L_p^w(G)\circ X$ is dense in $X$. For this purpose, it is
sufficient to show that closure of $Y$ contains indicator functions of all
compact sets. Let $A\subset G$ be a compact set. Take a relatively compact
neighborhood of identity $V=V^{-1}$, then $I_V^\trg = I_{V^{-1}} \in \L_p^w(G)$,
$I_{AV}\in X$ because both $V$ and $AV$ are relatively compact.
Put $f_V=I_{AV}*I_V/\mu(V)=I_V^\trg\circ I_{AV}/\mu(V)$, $f_V\in Y$.
Easy to check that
$$
I_A\le f_V \le I_{AV^2},
$$
whence $\|I_A-f_V\| \le \|I_A-I_{AV^2}\|$. At the same time
$$
\|I_A-I_{AV^2}\|_\pw^p = \int_{AV^2\setminus A} w^p(t)dt,
$$
what tends to zero as $\mu(V)\to0$ due to the fact that $w^p$ is locally summable.
This means that $I_A$ belongs to the closure of $Y$, and it follows that $Y=X$.

Thus, for every $\phi\in Y$ there are $f\in \L_q^{w^{-1}}(G)$, $g\in \L_p^w(G)$
such that $\phi=g\circ f=f*g^\trg$. But now
\begin{gather*}
|\phi(x)| = |(f*g^\trg) (x)| = \big|\int f(t)g(x^{-\!1}t)dt\big|
\le\\
\le \|f\|_{q,w^{-\!1}}\|\trans {x^{-\!1}}{\!g}\|_{p,w}
 \le L_x \|f\|_{q,w^{-\!1}}\|g\|_{p,w}.
\end{gather*}
This means that $\phi/L \in \L_\infty(G)$. As $Y=\L_q^{w^{-1}}(G)=w\L_q(G)$,
we get for all $\psi\in \L_q(G)$ that $\psi w/L\in \L_\infty$.
Since on every compact set $F$ the function $w$ is bounded away from zero
and $L$ is bounded, we get that $\L_q(F)\subset \L_\infty(F)$,
what is possible for finite $F$ only. This contradicts
assumption that $G$ is not discrete.
\epr

\end{document}